\documentclass[12pt]{amsart}
\usepackage{amsfonts,amsmath,amssymb,amsthm}
\usepackage{graphicx}
\usepackage{pgf,tikz}
\usetikzlibrary{arrows}
\usepackage[letterpaper, left=2.5cm, right=2.5cm, top=2.5cm,
bottom=2.5cm,dvips]{geometry}
\usepackage[T1]{fontenc}

\newtheorem{thrm}{Theorem}
\theoremstyle{plain}

\newtheorem{lem}{Lemma}
\newtheorem{cor}{Corollary}

\begin{document}
	\title[Laurent graph polynomials]{Laurent polynomials of planar graphs}
	
	\author{Mariusz Zaj\k{a}c}
	\address{Faculty of Mathematics and Information Science, Warsaw University
		of Technology, Koszykowa 75, 00-662 Warsaw, Poland}
	\email{M.Zajac@mini.pw.edu.pl}

	\begin{abstract}
	   By using Laurent graph polynomials instead of the usual ones, i.e. by allowing
	   negative powers of the variables, we simplify an existing method of determining
	   the Alon-Tarsi numbers of planar graphs.
	\end{abstract}
	
	\maketitle



The present note is closely related to the recent work on the Alon-Tarsi number
of planar graphs (\cite{Z} and \cite{GZ}),  
and uses some of the notation and terminology introduced in those papers
without further explanation.

Here are the main results: 
the graph polynomial 
of every planar graph $G$
$$
P_G=\prod_{xy\in E(G)}\left(x-y\right)
$$ 
has a nonvanishing monomial of degree at most
 $4$ with respect to every variable (Zhu, \cite{Z});  every planar graph $G$
contains such a matching $M$ that the polynomial
$P_{G-M}$
has a nonvani\-shing monomial of degree at most
$3$ with respect to every variable (Grytczuk and Zhu, \cite{GZ}).

\section*{A new approach to Zhu's theorem}

In order to carry out the inductive argument Zhu \cite{Z} formulates and proves,
following \cite{T}, a stronger fact:\\

(*) Let $G$ be a planar near-triangulation, and let $e = ab$ be a fixed edge in the outer cycle
of $G$. The polynomial $P_{G-e}$ contains then a nonvanishing 
monomial $N$ (i.e. one whose coefficient is a nonzero integer when computed 
over $\mathbb{Q}$, hence nonzero over an arbitrary field $\mathbb{F}$) with degrees:\\
(i) $\deg_a(N) = \deg_b(N) = 0$,\\
(ii) $\deg_v(N) \leqslant 2$ for every other $v$ in the outer cycle,\\
(iii) $\deg_u(N) \leqslant 4$ for all interior vertices $u$.\\

Let us introduce some modifications:
\begin{itemize}
	\item polynomials belonging to $\mathbb{F}[x,x^{-1},y,y^{-1},\ldots]$,
	usually called Laurent polynomials, will be considered instead of just $\mathbb{F}[x,y,\ldots]$;
	\item the \textbf{modified polynomial} of $G$ will be
	$$
	Q_G=\prod_{xy\in E(G)}\left(\frac1x-\frac1y\right);
	$$ 
	\item 
	a polynomial with no negative exponents will be called {\bf nice}; this yields
	a natural partial order: $M \succcurlyeq N$ 
	(`$M$ is at least as nice as $N$')
	whenever $\frac{M}{N}$ is nice (the symbol $\succcurlyeq$ will also denote 
	the relation  between the multidegrees, e.g.
	$(3,1,-1) \succcurlyeq (1,0,-1)$ will mean the same as
	$x^3yz^{-1} \succcurlyeq xz^{-1}$);
	\item
	a polynomial is {\bf good} if it has at least one nice monomial, 
	and {\bf bad} otherwise;
	\item 
	$P_1 \simeq P_2$ will mean that $P_1-P_2$ is bad; it is an equivalence relation
	because the set of bad polynomials is a linear space, and for the same reason no good 
	polynomial is in relation $\simeq$ with a bad one;

\end{itemize}

	Monomials in the graph polynomials will be further compared with \textbf{reference monomials}, e.g. 
for the graph $G$ discussed in (*) we will write $$R_{G-e}=v_1^2v_2^2\ldots v_k^2 u_1^4u_2^4\ldots u_l^4$$ and formulate the statement as follows: the polynomial $Q_{G-e}$ contains 
a nonvanishing monomial at least as nice as  $R_{G-e}^{-1}$, or in other words:

\begin{thrm} With all the above notations and assumptions \\
(*)	the polynomial $Z_{G-e} = Q_{G-e} \cdot R_{G-e}$ is good.
\end{thrm}

\noindent
\textbf{Proof} (induction on  $|V(G)|$). If $V(G) = \{a,b,c\}$, then
$$
Z_{G-e} = Q_{G-e} \cdot R_{G-e} = \left(\frac1a-\frac1c\right)\left(\frac1b-\frac1c\right)
\cdot c^2= 1 -a^{-1}c - b^{-1}c + a^{-1}b^{-1}c^2,
$$
which is clearly good, and we can proceed to the induction step. 

If the outer cycle of $G$ has a chord $f=cd$, dividing $G$ into $G_1$ and $G_2$, where
w.l.o.g. $e \in E(G_1)$ and $f \in E(G_1) \cap E(G_2)$, then
$$Z_{G-e} = Z_{G_1-e} \cdot Z_{G_2-f}$$
By the inductive assumption $Z_{G_1-e}$ has a nice monomial $M_1$ without variables
$a$ and $b$, while $Z_{G_2-f}$
has a nice monomial $M_2$ without variables $c$ i $d$. We claim that the nice monomial $M_1 \cdot M_2$ without variables
$a$ and $b$
appears in the product $Z_{G-e}$.

Indeed, in principle  $M_1 \cdot M_2$ could reduce with some other monomial product
(the product of two good polynomials need not be good, e.g.
$(1+xy^{-1})(1-x^{-1}y)=$ \linebreak$=xy^{-1}-x^{-1}y$), but that cannot happen here,
for if $N_1 \cdot N_2=-M_1 \cdot M_2$ then no variable can appear in 
$\frac{N_1}{M_1}=-\frac{M_2}{N_2}$ with a negative exponent:\\
the vertices not in $V(G_1)$ cannot, because they do not at all appear in $\frac{N_1}{M_1}$,\\
the vertices not in $V(G_2)$ cannot, because they do not at all appear in $\frac{M_2}{N_2}$,\\
and finally $c$ i $d$ cannot, because their exponent in $M_2$ is $0$, and
in $N_2$ it is at most 0 by definition of the modified graph polynomial.\\
A contradiction results, because the total degree of $\frac{N_1}{M_1}$ must be $0$,
since the graph polynomials are always homogeneous.

Now assume that the outer cycle of $G$ has no chord. Let $v$ and $t$ be the two 
consecutive vertices on the outer cycle following $a$ and $b$ (possibly $t=a$)
and let the neighbors of $v$ be called $b,x_1,\ldots,x_k,t$. Denote also $G'=G-v$.

Now we have
$$
Q_{G-e} = Q_{G'-e} \cdot \left(\frac1b-\frac1v\right)
\left(\frac1{x_1}-\frac1v\right) \ldots \left(\frac1{x_k}-\frac1v\right) 
\left(\frac1t-\frac1v\right),
$$
and
$$
R_{G-e} = R_{G'-e} \cdot x_1^2 \ldots x_k^2 v^2,
$$
because the allowed degree grows from 2 to 4 at $x_1,\ldots,x_k$,
0 ($v \not\in V(G')$) to 2 at $v$, and elsewhere remains unchanged.

Hence
$$
Z_{G-e} = Z_{G'-e} \cdot \left(\frac1b-\frac1v\right)
\left(\frac1{x_1}-\frac1v\right) \ldots \left(\frac1{x_k}-\frac1v\right) 
\left(\frac1t-\frac1v\right)x_1^2 \ldots x_k^2 v^2.
$$
In the next step let us notice that after multiplying out the above
expression we can recognize as bad and discard all monomials for which 
a negative power of either $b$ or $v$ comes from the factor to the right of 
the multiplication dot
(in $Z_{G'-e}$ there is no variable $v$, and $b$ may only appear in
a negative power, as there is no $b$ in $R_{G'-e}$). Thus
$$
Z_{G-e} \simeq - Z_{G'-e} \cdot \frac1v
\left(\frac1{x_1}-\frac1v\right) \ldots \left(\frac1{x_k}-\frac1v\right) 
\left(\frac1t-\frac1v\right)x_1^2 \ldots x_k^2 v^2 = 
$$
$$
= - Z_{G'-e} \cdot \frac{(v-x_1)\ldots(v-x_k)(v-t)}{v^k}
\cdot \frac{x_1 \ldots x_k}t= 
$$
$$
= - Z_{G'-e} \cdot \frac{v^{k+1}-(x_1+\ldots +x_k+t)v^k+ \ldots }{v^k}
\cdot \frac{x_1 \ldots x_k}t \simeq 
$$
$$
\simeq Z_{G'-e} \cdot (t+ x_1+\ldots +x_k - v)
\cdot \frac{x_1 \ldots x_k}t.
$$

The theorem is now a direct consequence of the following

\begin{lem}
	Let the polynomial $U(t,x_1,\ldots, x_k)$ contain a nonvanishing monomial $M$ of degree 
	$(\alpha,\beta_1,\ldots, \beta_k)$. Then the polynomial 
	$U \cdot (t+x_1+\ldots+ x_k-v)$ contains a nonvanishing monomial of degree
	$\succcurlyeq (\alpha+1,\beta_1-1,\ldots, \beta_k-1,0)$ (consequently,
	$$
	V(t,x_1,\ldots, x_k,v)=U(t,x_1,\ldots, x_k) \cdot (t+x_1+\ldots+ x_k-v)
	\cdot \frac{x_1 \ldots x_k}t
	$$
	contains a nonvanishing monomial of degree $\succcurlyeq (\alpha,\beta_1,\ldots, \beta_k,0)$).
\end{lem}

\noindent
\textbf{Proof.} The obvious candidate for the required monomial is $M \cdot t$,
and if $M \cdot t$ happens to reduce with some $N \cdot x_i$, then the assertion is satisfied 
by $N \cdot v$, which is of the same degree as $M \cdot tx_i^{-1}v$,
and cannot reduce with anything because the variable $v$ does not appear
anywhere else. \hfill$\square$\\

A direct application of Lemma 1 ends the inductive step: if $Z_{G'-e}$ has a nice
monomial, then $Z_{G-e}$ has a monomial that is at least as nice. \hfill$\square$ \\

Let us marginally remark that although
the case when the outer cycle of $G$ is a triangle
need not be considered separately, showing that $Z_{G-e}$ is good becomes then trivial as also
$t=a$ can be excluded from the denominator:
$$
Z_{G-e} \simeq Z_{G'-e} \cdot (t+ x_1+\ldots +x_k - v)
\cdot \frac{x_1 \ldots x_k}t  \simeq Z_{G'-e} \cdot x_1 \ldots x_k.
$$

\pagebreak

\section*{A stronger result by Grytczuk and Zhu}

Grytczuk and Zhu suggested in \cite{GZ} a strengthening of (*), which
can be rephrased as follows:  \\

\begin{thrm} 
If $G$ is a planar near-triangulation, and if $e = ab$ is a fixed edge in the outer cycle
of $G$, then there exists such a (possibly empty) matching $\mathcal{M}$ in $G$
consisting of the edges  $r_ig_i \in E(G)$, that \\
(a) red vertices ($r_i$), green vertices ($g_i$), $a$ and $b$ are pairwise distinct;\\
(b) no red vertex lies on the outer cycle of $G$; \\
(**) the polynomial $P_{G-e-\mathcal{M}}$ contains a nonvanishing 
(over an arbitrary field $\mathbb{F}$)
monomial $M$  with degrees:\\
(i) $\deg_a(M) = \deg_b(M) = 0$,\\
(ii) $\deg_g(M) \leqslant 1$ for every green vertex $g$ of the outer cycle,\\
(iii) $\deg_v(M) \leqslant 2$ for every other vertex $v$ of the outer cycle,\\
(iv) $\deg_u(M) \leqslant 3$ for every interior vertex $u$;\\
(***) the polynomial $P_{G-e}$ contains a nonvanishing 
(over an arbitrary field $\mathbb{F}$)
monomial $N$  with degrees:\\
(i) $\deg_a(N) = \deg_b(N) = 0$,\\
(ii) $\deg_g(N) \leqslant 1$ for every green vertex $g$ of the outer cycle,\\
(iii) $\deg_v(N) \leqslant 2$ for every other vertex $v$ of the outer cycle,\\
(iv)  $\deg_u(N) \leqslant 3$ for every non-red interior vertex $u$,\\
(v) $\deg_r(N) \leqslant 4$ for every red interior vertex $r$.\\
\end{thrm}

The proof of (**)  is analogous to that of (*), but it requires a minor change in the notation: the neighbors of $v$ other than $b$ and $t$ will be called $x_1, \ldots, x_k$
(green) and $y_1,\ldots, y_l$ (non-green, none of them is red). 
Now if we do not change the matching on adding $v$ to $G'$ 
(i.e. $\mathcal{M}=\mathcal{M}'$), then
$$
R_{G-e-\mathcal{M}} = R_{G'-e-\mathcal{M}'} \cdot x_1^2 \ldots x_k^2 y_1 \ldots y_l v^2,
$$
because at $x_1, \ldots, x_k$ the allowed degree grows from 1 to 3, 
at $y_1,\ldots, y_l$ from 2 to 3, and at $v$ from 0 (no vertex in $G'$) to 2. 
Just as in the proof of Theorem 1
\begin{equation}\label{c1}
Z_{G-e-\mathcal{M}} \simeq Z_{G'-e-\mathcal{M}'} \cdot (t+ x_1+\ldots +x_k + y_1+\ldots +y_l - v)
\cdot \frac{x_1 \ldots x_k}t.
\end{equation}

If we augment the matching by putting $\mathcal{M}=\mathcal{M}' \cup \{vy_1\}$, ($v$ will be green, $y_1$ will be red), then 
$$
R_{G-e-\mathcal{M}} = R_{G'-e-\mathcal{M}'} \cdot x_1^2 \ldots x_k^2 y_1 \ldots y_l v,
$$
$$
Q_{G-e-\mathcal{M}} =
$$ $$ 
Q_{G'-e-\mathcal{M}'} \cdot \left(\frac1b-\frac1v\right)
\left(\frac1{x_1}-\frac1v\right) \ldots \left(\frac1{x_k}-\frac1v\right)
\left(\frac1{y_2}-\frac1v\right) \ldots \left(\frac1{y_l}-\frac1v\right) 
\left(\frac1t-\frac1v\right).
$$
Discarding again the terms with $b$ and $v$ in the denominator, we obtain:

\begin{equation}\label{c2}
Z_{G-e-\mathcal{M}} \simeq - Z_{G'-e-\mathcal{M}'} \cdot y_1 \cdot \frac{x_1 \ldots x_k}t.
\end{equation}

The sum of the right sides of the equation (1) and $l$ equations (2) for $y_1,\ldots,y_l$
equals
$$
Z_{G'-e-\mathcal{M}'} \cdot (t+ x_1+\ldots +x_k - v)
\cdot \frac{x_1 \ldots x_k}t
$$
and is a good polynomial by Lemma 1, so at least one of those right sides
is good, which ends the proof of (**). \\

We will show (***) in the same way: if we do not change the matching on adding $v$ to $G'$ ($\mathcal{M}=\mathcal{M}'$), then, as above,
$$
R_{G-e} = R_{G'-e} \cdot x_1^2 \ldots x_k^2 y_1 \ldots y_l v^2,
$$
and if we put $\mathcal{M}=\mathcal{M}' \cup \{vy_1\}$, ($v$ is green, $y_1$ is red), then 
$$
R_{G-e} = R_{G'-e} \cdot x_1^2 \ldots x_k^2 y_1^2 y_2 \ldots y_l v,
$$
because at $x_1, \ldots, x_k$ the allowed degree grows from 1 to 3, at $y_1$ from 2 to 4,
at $y_2,\ldots, y_l$ from 2 to 3, and finally at $v$ from 0 to 1. 

In the former case, as before,
\begin{equation}\label{c3}
Z_{G-e} \simeq Z_{G'-e} \cdot (t+ x_1+\ldots +x_k + y_1+\ldots +y_l - v)
\cdot \frac{x_1 \ldots x_k}t,
\end{equation}
while in the latter
$$
Q_{G-e} =
$$ $$ 
Q_{G'-e} \cdot  \left(\frac1b-\frac1v\right)
\left(\frac1{x_1}-\frac1v\right) \ldots \left(\frac1{x_k}-\frac1v\right)
\left(\frac1{y_1}-\frac1v\right) \ldots \left(\frac1{y_l}-\frac1v\right) 
\left(\frac1t-\frac1v\right),
$$
and
\begin{equation}\label{c4}
Z_{G-e} \simeq - Z_{G'-e} \cdot y_1 \cdot \frac{x_1 \ldots x_k}t
\end{equation}
and we conclude the proof as in (**). We also see that the same edge $vy_i$
can be used for enlarging the matching in (**) and (***), because the formulas 
(3) and (4) are the same as (1) i (2). Consequently, at the end of the procedure
the matchings in cases (**) and (***) can be assumed to be equal.
\hfill$\square$ \\


\section*{A remark on $K_5$-minor-free graphs}

The following fact is a direct consequence of Theorem 2:

\begin{thrm} 
	If $G$ is a planar triangulation with the outer cycle $\Delta$,
	then there exists such a (possibly empty) matching $\mathcal{M}$ in $G$
	consisting of the edges  $r_ig_i \in E(G)$, that
	no red ($r_i$) or green ($g_i$) vertices belong to $V(\Delta)$, and\\
	(**) the polynomial $P_{G-E(\Delta)-\mathcal{M}}$ contains a nonvanishing 
	(over an arbitrary field $\mathbb{F}$)
	monomial $M$  with degrees:\\
	(i) $\deg_v(M) = 0$ for every vertex $v \in V(\Delta)$,\\
	(ii) $\deg_u(M) \leqslant 3$ for every vertex $u \not\in V(\Delta)$;\\
	(***) the polynomial $P_{G-E(\Delta)}$ contains a nonvanishing 
	(over an arbitrary field $\mathbb{F}$)
	monomial $N$  with degrees:\\
	(i) $\deg_v(N) = 0$ for every vertex $v \in V(\Delta)$,\\
	(ii)  $\deg_u(N) \leqslant 3$ for every non-red vertex $u \not\in V(\Delta)$,\\
	(iii) $\deg_r(N) \leqslant 4$ for every red vertex $r$.\\
\end{thrm} 

For a proof assume that $V(\Delta)=\{a,b,c\}$ and apply Theorem 2 to the graph
$G'=G-c$. \hfill$\square$ \\

\pagebreak

The above fact is important because of the structural theorem describing 
the $K_5$-minor-free graphs (\cite{Wagner}), 
which states that every edge-maximal graph with no $K_5$ minor
can be obtained from planar triangulations and one specific non-planar graph $V_8$
($V_8$ is a 3-regular graph on 8 vertices, obtained from $C_8$ by joining
pairs of opposite vertices) by glueing along cliques of at most 3 vertices.

From Theorems 2 and 3, as well as the obvious fact that in $P_{V_8}$ all exponents
equal at most 3, we can therefore deduce: 

\begin{thrm} 
	If $G$ is a $K_5$-minor-free graph,
	then there exists such a (possibly empty) matching $\mathcal{M}$ in $G$
	consisting of the edges  $r_ig_i \in E(G)$, that\\
	(**) the polynomial $P_{G-\mathcal{M}}$ contains a nonvanishing 
	(over an arbitrary field $\mathbb{F}$)
	monomial $M$  with degree $\deg_u(M) \leqslant 3$ for every vertex $u \in V(G)$;\\
	(***) the polynomial $P_{G}$ contains a nonvanishing 
	(over an arbitrary field $\mathbb{F}$)
	monomial $N$  with degrees:\\
	(i)  $\deg_u(N) \leqslant 3$ for every vertex $u \in V(G)$ that is not red,\\
	(ii) $\deg_r(N) \leqslant 4$ for every red vertex $r \in V(G)$.\\
\end{thrm} 

Applying Alon's celebrated Combinatorial Nullstellensatz (see \cite{AlonCN}), we easily obtain
the following corollaries:

\begin{cor}
	If $G$ is a $K_5$-minor-free graph,
	then there exists such a  
	matching $\mathcal{M}$ in $G$,
	that $G-\mathcal{M}$ is $4$-choosable.
\end{cor}

\begin{cor}
	If $G$ is a $K_5$-minor-free graph on $n$ vertices,
	then there exists such a subset $A \subset V(G)$
	with $|A|<\frac n2,$ that $G$ is list colorable
	assuming that at all vertices there is a list of $4$ 
	or $5$ colors, with $5$ colors appearing only
	at the vertices that belong to $A$.
\end{cor}

\end{document}